\documentclass[12pt]{article}
\usepackage{amsfonts,amssymb,latexsym,epsfig,amsmath}
\usepackage[english]{babel}
\usepackage{amssymb}
\usepackage{amsmath}

\newtheorem{theor}{Theorem}[section]
\newtheorem{lemma}{Lemma}[section]
\newtheorem{cor}{Corollary}[section]

\newtheorem{remark}{Remark}[section]

  \newcommand{\sn}{~\hfill\hbox{$\square$}}
\newcommand{\n}{\hbox{$\scriptstyle\hbox{\rm
     l\kern-0.22em N}$}}
\newcommand{\N}{\hbox{$ I\!\!N$}}

\renewcommand{\r}{{\mathbb R}}
\newcommand{\R}{\r^{n}}

\renewcommand{\O}{{\mathcal{O}}}
\newcommand{\e}{\varepsilon}

\newcommand{\rec}[1]{{(\ref{#1})}}

\def\PS #1 #2{\langle #1, #2 \rangle}
\def\YtoXandEPStoZERO{{\displaystyle{\mathop{\scriptstyle{y \to
x}}_{\varepsilon\to 0}}}}

\def\1{1\!\!1}

\def\Ob{{\overline{\mathcal{O}}}}
\def\dO{{\partial{\mathcal{O}}}}
\def\Kb{{\overline K}}
\def\ou{\overline u}
\def\uu{\underline u}

\def\ow{\overline w}

\def\eps{\varepsilon}
\def\control{(\alpha_t)_t}
\def\Ex{I\!\!E_x}

 1
\begin{document}

\title{\bf On the Boundary Ergodic Problem for Fully Nonlinear Equations in Bounded Domains with General Nonlinear Neumann Boundary Conditions}

\author{Guy\ Barles$^{(1)}$ \&
Francesca Da Lio$^{(2)}$}

\addtocounter{footnote}{1} \footnotetext{Laboratoire de
Math\'ematiques et Physique Th\'eorique.  Universit\'e de Tours. Facult\'e des Sciences et Techniques, Parc de Grandmont, 37200
Tours, France.} \addtocounter{footnote}{1}

\footnotetext{Dipartimento di Matematica.  Universit\`a di Torino. Via Carlo Alberto 10,  10123 Torino, Italy.}


 \maketitle

\begin{abstract} We study nonlinear Neumann type boundary value problems related to ergodic phenomenas. The particularity of these problems is that the ergodic constant appears in the (possibly nonlinear) Neumann boundary conditions. We provide, for bounded domains, several results on the existence, uniqueness and properties of this ergodic constant.

\end{abstract}

\markboth{\textsc{Barles \& Da Lio}}{\textsc{Boundary Ergodic Problems} }

\section{Introduction}

In this article, we are interested in what can be called ``boundary ergodic control problems'' which lead us to solve the following type of fully nonlinear elliptic equations associated with nonlinear Neumann boundary conditions
\begin{eqnarray}
F(x,Du,D^2u)=\lambda & \mbox{in ${\mathcal{O}}$,} \label{eqn} \\
L(x,Du)=\mu & \mbox{on $\partial{\mathcal{O}} $,} \label {bc}
\end{eqnarray}
where, say, $\O\subset\R$ is a smooth domain, $F $ and $L$ are, at least, continuous functions defined respectively on $\overline{\mathcal{O}} \times\R\times {\mathcal{S}}^n$ and $\overline{\mathcal{O}} \times\R$ with values in $\r$, where ${\mathcal{S}}^n$ denotes the space of real, $n\times n$, symmetric matrices. More precise assumptions on $F$ and $L$ are given later on.

The solution $u$ of this nonlinear problem is scalar and $Du$, $D^2u$ denote respectively  gradient and Hessian matrix of $u.$
Finally, $\lambda$, $\mu$ are constants~: $\mu$, which is called below the ``boundary ergodic cost'', is part of the unknowns while $\lambda$ is mainly here considered as a given constant for reasons explained below.

In order to justify the study of such problems, we first concentrate only on the  equation (\ref{eqn}), without boundary condition, i.e. on the case when $\O=\R$. In this framework, under suitable assumptions on $F$, the typical result that one expects is the following~: there exists a unique constant $\lambda$ such that (\ref{eqn}) has a bounded solution. Such results were first proved for first-order equations by Lions, Papanicolaou \& Varadhan \cite{LPV} (see also Concordel \cite{Conc1}) in the case of periodic equations and solutions. Recently, Ishii \cite{I1} generalizes these results in the almost periodic case. General results for second-order equations in the periodic setting are proved by Evans \cite{E1,E2}. Results in the evolution case, when the equation is periodic both in space and time, were also obtained recently by Souganidis and the first author \cite{BS2}~: the methods of \cite{BS2}, translated properly to the stationary case, are the one who would lead to the most general results in the case of second-order equations. All these results which hold for general equations without taking advantage of their particularities, are complemented by more particular results in the applications we describe now.

The first application concerns the so-called ergodic control
problems (either in the deterministic or stochastic case). We
refer to Bensoussan \cite{B} for an introduction to such problems
and to Bensoussan \& Frehse \cite{BFas}, Bagagiolo, Bardi \&
Capuzzo Dolcetta \cite{BBC},  Arisawa \cite{A1,A2}, Arisawa \&
Lions \cite{AL} for further developments in the $\R$ case and with
different types of pde approaches. In this framework, (\ref{eqn})
is the Bellman Equation of the ergodic control problem, $\lambda$
is the ergodic cost and the solution $u$ is the value function of
the control problem. In this case, both the uniqueness of
$\lambda$ and of $u$ -- which can hold only up to an additive
constant -- is interesting for the applications. But  it is rather
easy to obtain the uniqueness of $\lambda$ in general, while the
uniqueness of $u$ can be proved only in the uniformly elliptic
case and is generally false.

A second motivation to look at such problems is the asymptotic behavior as $t \to \infty$ of solutions of the evolution equation
\begin{equation}
\label{evol}
u_t + F(x,Du,D^2 u) = 0 \quad\hbox{in  }\R \times (0,+\infty)\; .
\end{equation}

A typical result here is the following~: if there exists a unique $\lambda$ such that (\ref{eqn}) has a solution (typically in the bounded solutions framework), then one should have
$$ \frac{u(x,t)}{t}Ê\to \lambda \quad\hbox{locally uniformly as  }Êt \to \infty\; .$$
Therefore the ergodic constant governs the asymptotic behavior of the associated evolution equation and in good cases, one can even show that
$$ u(x,t) -  \lambda t \to u_\infty (x)   \quad\hbox{as  }Êt \to \infty\; ,$$
where $u_\infty$ solves (\ref{eqn}).

Such results were obtained recently, for first-order equations, by Fathi \cite{Fa1, Fa2, Fa3} and Namah \& Roquejoffre \cite{NR} in the case when $F$ is convex in $Du$ ; these results were generalized and extended to a non-convex framework in Barles \& Souganidis \cite{BS1}. To the best of our knowledge, there is not a lot of general results in the case of second-order equations~: the uniformly elliptic case seems the only one which is duable through the use of the Strong Maximum Principle and the methods of \cite{BS2} which are used in the paper to prove the convergence to space-time periodic solutions but which can be used to show the convergence to solutions of the stationary equations.

The third and last application (and maybe the most interesting one) concerns homogenization of elliptic and parabolic pdes. This was the motivation of Lions, Papanicolaou \& Varadhan \cite{LPV} to study these types of ergodic problems as it was also the one of Evans \cite{E1,E2}. The ergodic problem is nothing but the so-called ``cell problem'' in homogenization theory, $\lambda$ being connected to the ``effective equation''. We also refer the reader to Concordel \cite{Conc2}, Evans \& Gomes \cite{EG}, Ishii \cite{I1} for results in this direction. The connections between ergodic problems and homogenization are studied in a systematic way in Alvarez \& Bardi \cite{AB1,AB2} and completely clarified.

Of course, the same questions have been studied in bounded (or unbounded) domains with suitable boundary conditions. For first-order equations, Lions \cite{L} studies the ergodic problem in the case of homogeneous Neumann boundary conditions, while Capuzzo Dolcetta \& Lions \cite{icdpll} study it in the case of state-constraints boundary conditions. For second-order equations, we refer the reader to Bensoussan \& Frehse \cite{BFnbc} in the case of homogeneous Neumann boundary conditions and to Lasry \& Lions \cite{LaLi}  for state-constraints boundary conditions. It is worth pointing out that in all these works, the constant $\mu$ does not appear and the authors are interested in the constant $\lambda$ instead.

The first and, to the best of our knowledge, only work where the
problem of the constant $\mu$ appears, is the one of Arisawa
\cite{A3}. In this work, she studies two different cases~: the
case of bounded domains { which we consider here and the case of
half-space type domains which contains different difficulties; we
address this problem in a forthcoming work in collaboration with
P.L. Lions and P.E. Souganidis.Ê In the case of bounded domains,
we improve her results in several directions~: generality and
regularity of the equation and boundary condition, possibility of
obtaining results in degenerate cases, uniqueness in more general
frameworks, interpretation in terms of stochastic control problems
and connections of these types of ergodic problems with large time
behavior of solutions of initial value problem with Neumann
boundary conditions. We are able to do so since we use softer
viscosity solutions' methods.

It is worth pointing out that the role of the two constants are
different~: our main results say that, for any $\lambda$, there
exists a unique constant $\mu:= \mu(\lambda)$ for which
(\ref{eqn})-(\ref{bc}) has a bounded solution. Therefore the role
played previously by $\lambda$ is now played by $\mu$. To prove
such a result, we have to require some uniform ellipticity
assumption on $F$, not only in order to obtain the key estimates
which are needed to prove the existence of the solution $u$ but
also because $\mu$ can play its role only if the boundary
condition is ``seen in a right way by the equation". Indeed, the
counter-example of Arisawa \cite{A3}, p. 312, shows that otherwise
$\mu$ cannot be unique. This vague statement is partly justified
in Section~\ref{escp}.

The proof of the existence of the solution relies on the
$C^{0,\alpha}$ estimates proved in \cite{BDL};  in order to have
an as self-contained paper as possible, we describe these results
in the Appendix. Here also the uniform ellipticity of $F$ plays a
role, at least in the case when the Neumann boundary condition is
indeed nonlinear. But if the boundary condition is linear, some
less restrictive ellipticity assumptions on $F$ can be made~: this
is the reason why we distinguish the two cases below.

An other question we address in this paper, are the connections with the large time behavior of the solutions of the  two different type of evolution problems
\begin{eqnarray}
 v_t +  F(x,Dv,D^2 v)& =& \lambda \quad \hbox{in $\mathcal{O} \times (0,+\infty)$,} \label{ev-eqn1}\\
   L(x,Dv) & = &  \mu  \quad\hbox{on $\dO \times (0,+\infty) $.}                 \label {ev-bc1}
\end{eqnarray}
and
\begin{eqnarray}
 w_t +  F(x,Dw,D^2 w)& =& 0 \quad \hbox{in $\mathcal{O} \times (0,+\infty)$,} \label{ev-eqn2}\\
w_t +   L(x,Dw) & = & 0  \quad\hbox{on $\dO \times (0,+\infty) $.}                 \label {ev-bc2}
\end{eqnarray}

In the case of (\ref{ev-eqn1})-(\ref{ev-bc1}), we show that the ergodic constant $\mu(\lambda)$ is characterized as the only constant $\mu$ for which the solution $v$ remains bounded. In the case of (\ref{ev-eqn2})-(\ref{ev-bc2}), the expected behavior is to have $\displaystyle t^{-1}w(x,t)$ converging to a constant $\tilde \lambda$ which has to be such that (\ref{eqn})-(\ref{bc}) has a solution for $\tilde \lambda=\lambda = \mu(\lambda)$. We prove that, under suitable conditions, such a constant $\tilde \lambda$, i.e. a fixed point of the map $\lambda \mapsto \mu(\lambda)$, does exist and that we have the expected behavior at infinity for $w$.

Finally we consider the case when the equation is the Hamilton-Jacobi-Bellman Equation of a  stochastic control problem with reflection~: this gives us the opportunity to revisit the results on the uniqueness of $\mu$ in a \emph{degenerate} context and to provide a formula of representation for $\mu$.

The paper is organized as follows. In Section~\ref{S:nlbc}, we prove the existence of $u$ and $\mu$ in the case of nonlinear boundary conditions while in Section~\ref{S:lbc} we treat the linear case. In Section~ \ref{unimu}, we examine the uniqueness properties for $\mu$ together with its dependence in $\lambda$, $F$ and $L$; among the results of this part, there is the existence of $\tilde \lambda$. Section~\ref{ab} is devoted to present the results connecting the ergodic problem with the asymptotic behavior of solution of some nonlinear problem. Finally we study the connections with stochastic control problem with reflection in Section~\ref{escp}.

\section{The case of nonlinear boundary conditions}\label{S:nlbc}

To state our result, we use the following assumptions

\smallskip

\noindent{\bf (O1)} ${\mathcal{O}}$ is a bounded domain with a $W^{3,\infty}$ boundary.

\smallskip
We point out that such assumption on the regularity of the boundary is needed both in order to use the comparison results of \cite{b2} (here the $W^{3,\infty}$ regularity is needed) and the local $C^{0, \alpha}$-estimates of \cite{BDL} (here a $C^{2}$ regularity would be enough). 

We denote by $d$ the sign-distance function to $\partial{\mathcal{O}}$ which is positive in ${\mathcal{O}}$ and negative in $\R \setminus \overline{\mathcal{O}}$. If $x\in \partial{\mathcal{O}}$, we recall that $Dd(x)=-n(x)$ where $n(x)$ is the outward unit normal vector to $\partial{\mathcal{O}}$ at $x$. The main consequence of {\bf (O1)} is that $d$ is $W^{3,\infty}$ in a neighborhood of $\partial{\mathcal{O}}$.

Next we present the assumptions on $F$ and $L$.

\smallskip

\noindent {\bf (F1) (Regularity)} The function $F$ is locally Lipschitz continuous on $\Ob \times \R \times {\mathcal{S}}^n$ and there exists a constant $K>0$ such that, for any $x,y\in \Ob$, $p,q \in \R$, $M,N \in {\mathcal{S}}^n$
$$ |F(x,p,M) - F(y,q,N)|\leq K\left\{ |x-y|(1+|p|+|q| + ||M||+||N||) +|p-q| + ||M - N||\right\}\;. $$

\noindent {\bf (F2) (Uniform ellipticity)} There exists $\kappa >0$ such that, for any $x\in \Ob$, $p\in \R$, $M,N \in {\mathcal{S}}^n$ with $N\geq 0$
$$ F(x,p,M+N) - F(x,p,M) \leq -\kappa{\rm Tr}(N)\; .$$

\noindent {\bf (F3) } There exists a continuous function $F_\infty$ such that
$$t^{-1} F(x, tp,tM)\to F_\infty (x,p,M)\quad \hbox{locally uniformly, as  } t\to + \infty\; .$$

For the boundary condition $L$, we use the following assumptions.

\par

\medskip

\noindent {\bf (L1)} There exists $\nu>0$ such that, for all $(x,p)\in
\partial{\mathcal{O}}\times\r^n$ and $t >0$, we have
\begin{equation}\label{croiss}
L(x,p+t n(x))-L(x,p)\ge \nu t\,.
\end{equation}

\noindent {\bf(L2)} There is  a constant $\Kb >0$ such that, for all $x,y\in \partial{\mathcal{O}},$ $p,q\in\r^n$, we have
\begin{equation}
|L(x,p)-L(y,q)|\le \Kb \left[(1+ |p|+|q|)|x-y|+|p-q|
\right]\,.
\end{equation}

\noindent {\bf (L3)} There exists a continuous function $L_\infty$ such that
$$t^{-1} L(x, tp)\to L_\infty (x,p)\quad \hbox{locally uniformly, as  } t\to + \infty\; .$$

{ Before stating and proving the main result of this section, we
want to emphasize the fact that the above assumptions are very
well adapted for applications to stochastic control and
differential games: indeed {\bf (F1)-(L1)} are clearly satisfied
as soon as the dynamic has bounded and Lipschitz continuous drift,
diffusion matrix and direction of reflection and when the running
and boundary cost satifies analogous properties (maybe these
assumptions are not optimal but they are rather natural) while
{\bf (F3)-(L3)} are almost obviously satisfied because of the
structure of the Bellman or Isaac Equations (``sup'' or ``inf
sup'' of affine functions in $p$ and $M$)}.

Our result is the
\begin{theor}\label{ergnl} Assume   {\bf(O1)}, {\bf(F1)}-{\bf(F3)} and {\bf (L1)}-{\bf(L3)} then, for any $\lambda \inÊ\r$, there exists $\mu \in \r$ such that (\ref{eqn})-(\ref{bc}) has a continuous viscosity solution.
\end{theor}

\noindent{\bf Proof.} The proof follows the strategy of Arisawa \cite{A3}. For $0<\varepsilon \ll \alpha \ll 1$, we introduce the approximate problem
\begin{eqnarray}
F(x,D\tilde u,D^2\tilde u)+\e \tilde u=\lambda & \mbox{in ${\mathcal{O}}$,} \label{eqna} \\
L(x, D\tilde u) +\alpha\tilde u=0& \mbox{on $\partial{\mathcal{O}} $.} \label{bca}
\end{eqnarray}

\noindent{\bf 1.} It is more or less standard to prove that this problem has a unique continuous viscosity solution using the Perron's method of Ishii \cite{I2}
 and the comparison arguments of Barles \cite{b2}; the only slight difficulty comes from the $x$-dependence of $F$ which is a priori not sufficient to obtain a suitable comparison. In the Appendix, we explain why the usual approach does not work and we show  how to overcome this difficulty by borrowing ideas of Barles \& Ramaswamy \cite{GBMR}.

\noindent{\bf  2.} The next step consists in obtaining basic estimates on $\tilde u$. We drop the dependence of $\tilde u$ in $\e$ and $\alpha$ for the sake of simplicity of notations. To do so, we use the fact that ${\mathcal{O}}$ is bounded and therefore we can assume without loss of generality that ${\mathcal{O}} \subset \{x_1 >0\}$.

We introduce the smooth functions
$$ \ou (x) = C(2-\exp(-\gamma x_1))\quad , \quad  \underline {u} (x) = - C(2-\exp(-\gamma x_1))\; .$$
Notice that $ \uu < 0 <\ou$ on $\Ob$.

By using  {\bf (F1)} and  {\bf (F2)}, one sees that, for $\gamma$ and $C$ large enough, one has
\begin{eqnarray*}
F(x, D \ou , D^2\ou)&\ge& F(x,0,0)\\
&-&KC\gamma \exp(-\gamma x_1)+kC\gamma^2\exp(-\gamma x_1)>0\,,
\end{eqnarray*}
and
\begin{eqnarray*}
F(x, D \uu , D^2\uu)&\le& F(x,0,0)\\
&+&KC\gamma \exp(-\gamma x_1)-kC\gamma^2\exp(-\gamma x_1)<0\,.
\end{eqnarray*}

 Next we consider $\max_\Ob(\tilde u - \ou)$ and $\min_\Ob(\tilde u - \uu)$ which are achieved respectively at $x_1,x_2\in\Ob.$ Because of the above properties and since $\tilde u$ is a viscosity solution of (\ref{eqna})-(\ref{bca}),  these $\max$ and $\min$ cannot be achieved in ${\mathcal{O}}$ and, in any case, the ``$F$'' inequalities cannot hold. The ``$L$'' inequalities lead to the estimates
$$
\alpha \tilde u(x)\le \alpha(\overline u(x)-\overline u(x_1))+\sup_{\Ob}|L(x,D\overline u(x) )|\,,
$$
and
$$
\alpha \tilde u(x)\ge \alpha(\uu(x)-\uu(x_2))-\sup_{\Ob}|L(x,D\uu(x))|
$$
for every $x\in\Ob.$ Thus for some positive constant $\overline C(F,L)$ (depending on $F$ and $L$) we have
\begin{equation}\label{alphatildeu}
||\alpha \tilde u||_\infty \leq \overline C(F,L)\; .
\end{equation}
\noindent{\bf 3.} Let $x_0$ be any point of $\Ob$ and set $v(x) = \tilde u(x) - \tilde u(x_0)$ for $x \in \Ob$.
We claim that $v$ remains uniformly bounded as $\alpha$ tends to $0$ if $\e\ll\alpha$.
\par
To prove the claim, we argue by contradiction assuming that $M:=||v||_\infty \to \infty$ as $\alpha \to 0$ and we set $w(x) :=M^{-1}v(x)$. The function $w$ solves
\begin{eqnarray}
M^{-1} F(x,MDw,MD^2w)+\e w=M^{-1} \lambda - M^{-1} \e \tilde u(x_0) & \mbox{in ${\mathcal{O}}$,} \label{eqnb} \\
M^{-1} L(x, Dw) +\alpha w= -M^{-1} \alpha \tilde u(x_0) & \mbox{on $\partial{\mathcal{O}} $.} \label{bcb}
\end{eqnarray}
Moreover $||w||_\infty=1$ and $w(x_0)=0$.\par
Since $w$ is uniformly bounded, the $C^{0,\beta}$ regularity results and estimates of Barles \& Da Lio \cite{BDL} apply and therefore $w$ is uniformly bounded in $C^{0,\beta}$, for any $0<\beta<1$ (see also the Appendix, for a description of these results) .
\par
Using Ascoli's Theorem, one may assume without loss of generality that $w$ converges uniformly to some $C^{0,\beta}$-function $\ow$
and taking {\bf (F3)-(L3)} in account, the stability results for viscosity solutions implies that $\ow$ solves
\begin{eqnarray}
F_\infty (x,D\ow,D^2\ow)=0 & \mbox{in ${\mathcal{O}}$,} \label{eqnc} \\
L_\infty(x, D\ow) = 0 & \mbox{on $\partial{\mathcal{O}} $.} \label{bcc}
\end{eqnarray}
Moreover $||\ow||_\infty=1$ and $\ow(x_0)=0$.
\par
We are going to show now that all these properties lead to a contradiction by Strong Maximum Principle type arguments.  Since $\ow$ is continuous there exists $x\in \Ob$ such that $|\ow(x)|=1$.

We first remark that $F_\infty$ satisfies  {\bf (F1)-(F2)} as well
and is homogeneous of degree $1$; therefore the Strong Maximum
Principle of Bardi \& Da Lio \cite{BaDL} implies that necessarily
$x \in \dO$. In fact $-1 < \ow < 1$ in ${\mathcal{O}}$.
\par
We assume for example that $\ow(x)=1$, the other case being treated similarly. To conclude, we are going to use the following lemma.
\begin{lemma}\label{SMPbc} There exists $r>0$ and a smooth function $\varphi$ on $\overline B(x,r)$ such that $\varphi(x)=0$, $\varphi (y) > 0$ on $\dO \cap \overline B(x,r)\setminus\{x\}$
\begin{equation}
\label{Ffi}
F_\infty (y,D \varphi (y), D^2 \varphi(y)) >0 \quad\hbox{on  } \overline B(x,r)\; ,
\end{equation}
and
$$ D\varphi (x) = kn(x)\; ,$$
with $k>0$.
\end{lemma}

The proof of this lemma is given in the Appendix;  we show how to use it in order to conclude.
\par
Since $ D\varphi (x) = kn(x)$, we have $L_\infty(x,D\varphi (x))>0$. But $\varphi$ is smooth and therefore, by choosing $\theta < r$ small enough, we have also
\begin{equation}
\label{Lfi}
L_\infty (y,D \varphi (y)) >0 \quad\hbox{on  } \overline B(x,\theta)\cap \dO\; .
\end{equation}
On an other hand, by choosing
 $\tau>0$ small enough, we can have
 $ \ow(y)-\tau \varphi (y) < 1 = \ow(x)-\tau \varphi (x)$ for $y \in   \partial B(x,\theta)\cap \Ob$.
 Indeed, for $y$ close to $\dO$, $\varphi (y) >0$ while in  ${\mathcal{O}}$ we have  $\ow(y)<1$.
\par
We deduce from this property that, if we consider $\max_{ \overline B(x,\theta)\cap \Ob}( \ow-\tau \varphi)$, this maximum is
necessarely achieved in $ B(x,\theta)\cap \Ob$ and therefore it is a local maximum point of $ \ow-\tau \varphi$ but, taking in account the fact that $F_\infty$ and $L_\infty$ are homogeneous of degree 1, this is a contradiction with the inequalities  (\ref{Ffi})-(\ref{Lfi}).

\noindent{\bf  4.} From step~3, the functions $v$ are uniformly bounded and solve
\begin{eqnarray}
F(x,Dv,D^2v)+\e v= \lambda - \e \tilde u(x_0) & \mbox{in ${\mathcal{O}}$,} \label{eqnd} \\
L(x, Dv) +\alpha v= - \alpha \tilde u(x_0) & \mbox{on $\partial{\mathcal{O}} $.} \label{bcd}
\end{eqnarray}
Using again the regularity results of Barles \& Da Lio \cite{BDL} (see also the Appendix), we deduce that the functions $v$ are also uniformly bounded in $C^{0,\beta}$ for any $0<\beta<1$ and by Ascoli's Theorem, extracting if necessary a subsequence, we may assume that  they converge uniformly to a function $u \in C^{0,\beta}(\Ob)$. Moreover, since $\alpha \tilde u$ is also uniformly bounded, we can also extract a subsequence such that $- \alpha \tilde u(x_0)$ converges to some $\mu \in \r$.
\par
In order to conclude, we just pass to the limit in (\ref{eqnd})-(\ref{bcd}) with a choice of $\e$ such that $\e\alpha^{-1} \to 0$.\sn

\section{The case of linear boundary conditions}\label{S:lbc}

We consider in this section the case when $L$ is given by
\begin{equation}
\label{lbc} \langle Du,\gamma(x)\rangle  +g(x)=\mu \quad \mbox{on
$\partial{\mathcal{O}} $,}
\end{equation}
where the functions $\gamma$ and $g$ satisfies
\par
\medskip
\noindent {\bf (L1')} $g \in C^{0,\beta}(\dO)$ for some $0<\beta
\leq 1$ and $\gamma$ is a Lipschitz continuous function, taking
values in $\R$ and such that $\langle \gamma(x), n(x)\rangle
\ge\nu>0$ for any $x \in \dO$, where we recall that $n(x)$ denotes
the unit exterior normal vector to $\partial{\mathcal{O}}$ at $x$.
\par
\medskip
In this linear case, we are able to weaken the ellipticity assumption on $F$. In the following, for $q\in \R$, the notation $\hat q$ stands for $\displaystyle\frac{q}{|q|}$.

\par

\medskip

\noindent {\bf (F2') (partial uniform ellipticity)} There exists a Lipschitz continuous function $x \mapsto A(x)$, defined on $\Ob$ and taking value in the space of symmetric, definite positive matrix and $\kappa >0$ such that

(i)  for any $x\in \dO$, $A(x)\gamma(x) = n(x)$,

(ii)  for any $x\in \Ob$, $p\in \R\setminus\{0\}$, $M,N \in
{\mathcal{S}}^n$ with $N\geq 0$
$$ F(x,p,M+N) - F(x,p,M) \leq -\kappa\langle Nq,q\rangle + o(1)||N|| \; ,$$
with $q=\widehat {A^{-1}(x)p}$ and where $o(1)$ denotes a function of $|p|$ which converges to $0$ as $|p|$ tends to $+\infty$.
\par
\medskip
If $\gamma\equiv n$, this assumption is satisfied in particular if (formally)
$$F_M (x,p,M)  \leq -\kappa \hat p \otimes \hat p + o(1) \quad\hbox{a.e. in  }\Ob \times \R \times {\mathcal{S}}^n\,,$$
where, as above, $o(1)$ denotes a function of $|p|$ which converges to $0$ as $|p|$ tends to $+\infty$; this means a non-degeneracy property in the gradient direction, at least for large $|p|$.
 This corresponds to the choice $A(x)\equiv Id$. We recall that
 for all $p\in \R$, $p\otimes p$ denotes the symmetric matrix
 defined by $(p\otimes p)_{ij}=p_ip_j.$

In this case, unlike  the uniform elliptic case,  {\bf (F2')} is not enough to ensure a comparison property
for $F$, thus we add

\par

\medskip
\noindent {\bf (F4')} For any $\tilde K  >0$, there exists a
function $m_{\tilde K}: \r^+ \to \R$ such that $m _{\tilde K}
(t)\to 0$ when $t\to 0$ and such that, for all $\eta >0$

\noindent

\centerline{$ F(y,q,Y)-F(x,p,X)\leq m_{\tilde K} {\displaystyle
\left(\eta + |x-y|(1+\vert p \vert \vee \vert q\vert)
+\frac{|x-y|^2}{\varepsilon^2}\right)}$}

\medskip

\noindent

for all $x,y \in \Ob $, $p,q\in \R$ and for all  matrices $X, Y \in {\mathcal{S}}^n$ satisfying the
following properties
\begin{equation}\label{inegmatthbis}
-\frac{\tilde K}{\varepsilon^2 } Id
\leq \left(\begin{array}{cc} X & 0 \\ 0 & -Y
\end{array}\right) \leq \frac{\tilde K}{\varepsilon^2}
\left(\begin{array}{cc} Id
& -Id \\ -Id & Id
\end{array}\right) +  \tilde K\eta Id\, ,
\end{equation}
\begin{equation}
\vert p-q \vert \leq \tilde K \eta \varepsilon (1 + \vert p \vert \wedge
\vert q
\vert)\,,
\label{propgrad}
\end{equation}
\begin{equation}
\vert x-y \vert \leq \tilde K \eta \varepsilon .
\label{propxy}
\end{equation}

Our result is the

\begin{theor}\label{ergl} Assume  {\bf(F1)}-{\bf(F2')}-{\bf(F3)}-{\bf(F4')} and {\bf (L1')} then, for any $\lambda \inÊ\r$, there exists $\mu\in \r$ such that (\ref{eqn})-(\ref{bc}) has a continuous viscosity solution.
\end{theor}

We skip the proof since it follows readily the one of Theorem~\ref{ergnl}; we just point out that the key $C^{0,\beta}$-estimates follow from the linear case in \cite{BDL} (see also the Appendix) while the Strong Maximum Principle still holds under
{\bf(F2')} as we pointed it out in the Appendix.

\section{Uniqueness results for the boundary ergodic cost}\label{unimu}

In standard problems, the uniqueness of the ergodic cost is rather easy to obtain, while the uniqueness of the solution $u$ is a more difficult question. Here, even the uniqueness of $\mu$ is a non obvious fact because $\mu$ appears only in the boundary condition and clearly this boundary condition has to be sufficiently ``seen'' in order to have such a uniqueness property. The counter-example of Arisawa \cite{A3} for first-order equations shows that, in the cases where losses of boundary conditions can occur, $\mu$ is not unique in general.
\par
To state the uniqueness result, we introduce the following abstract assumption
\par
\medskip
\noindent{\bf (U1)} If $w$ is an upper semicontinuous viscosity subsolution of (\ref{eqn})-(\ref{bc}), there exists a sequence $(w_\e)_\e$ of upper semicontinuous functions such that ${\limsup}^*w_\e =w$ on $\Ob$,
 \footnote{We recall that the half-relaxed limit ${\limsup}^*w_\e$ is defined by~: $\displaystyle {\limsup}^*w_\e (x)=\limsup_{\YtoXandEPStoZERO}\,w_\e (y)$ for any $x \in \overline\O$.} satisfying in the viscosity sense
\begin{eqnarray}
F(x,Dw_\e,D^2 w_\e)\leq \lambda_\e < \lambda & \mbox{in ${\mathcal{O}}$,} \label{eqnw} \\
L(x,Dw_\e)\leq \mu + o_\e(1) & \mbox{on $\partial{\mathcal{O}} $.} \label {bcw}
\end{eqnarray}

Our result is the
\begin{theor}\label{unierg} Under the assumptions of either Theorem~\ref{ergnl} or \ref{ergl} and if {\bf (U1)} holds, if $u_1$ is a subsolution of (\ref{eqn})-(\ref{bc}) associated to $\lambda_1, \mu_1$ and if $u_2$ is a supersolution of (\ref{eqn})-(\ref{bc}) associated to $\lambda_2, \mu_2$ with $\lambda_1 \leq \lambda_2$ then necessarily $\mu_1 \geq \mu_2$.
In particular, for any $\lambda$, the boundary ergodic cost $\mu$ is unique.
\end{theor}
\noindent{\bf Proof.} We argue by contradiction assuming that $\mu_1 < \mu_2$.
\par
Let $u^\e_1$ be a continuous function associated to $u_1$ through assumption {\bf (U1)} with $\e$ choosen in such a way that
$$ L(x,Du^\e_1)\leq \mu \quad \mbox{on $\partial{\mathcal{O}} $,}$$
where $\mu := \displaystyle\frac{1}{2}(\mu_1 + \mu_2)$.
\par
We consider $\max_{\Ob \times \Ob}\,(u^\e_1(x)-u_2(y) -
\psi_\alpha (x,y))$ where for all $\alpha>0$  $\psi_\alpha $ is
the test-function built in \cite{b2} for the boundary condition
$L-\mu$ (we recall that this test-function depends only on the
boundary condition).
\par
Following readily the arguments of \cite{b2}, one is led to the inequalities
\begin{equation}\label{lambda1}
 F(x,p,X) \leq \lambda_{1,\e} < \lambda_1 \; ,
\end{equation}
\begin{equation}\label{lambda2}
 F(y,q,Y) \geq \lambda_2\; ,
\end{equation}
where $(p,X)\in D^{2,+}u^\e_1(\overline x)$ and $(q,Y)\in D^{2,-}u_2 (\overline y)$. Then either the standard comparison arguments or the arguments of \cite{GBMR} shows that
$$ F(x,p,X)-F(y,q,Y)\geq o_\alpha (1)\; ,$$
and hence, by subtracting the  inequalities \rec{lambda1}
and \rec{lambda2}, we have  $o_\alpha (1)\leq \lambda_{1,\e} - \lambda_2 < 0$. We get the contradiction by letting $\alpha$ tends to $0$. And the proof of the first part is complete.
\par
Of course, the uniqueness of the boundary ergodic cost follows since, if $u$ and $v$ are two solutions of (\ref{eqn})-(\ref{bc}) with the same $\lambda$ and roles with $\mu$, $\tilde \mu$ respectively, we can apply the above result with $u_1=u$, $\mu_1=\mu$, $\lambda_1=\lambda$ and $u_2=v$, $\mu_2=\tilde \mu$, $\lambda_2=\lambda$~: this yields $\mu_1 \geq \mu_2$. But using that the two solutions play symmetric roles, we deduce immediately $\mu=\tilde \mu$, i.e. the uniqueness of the ergodic cost.\sn

\begin{remark}\label{rem:unierg}{\rm
 As the proof shows it, the result ``$\mu_1 \leq \mu_2 \Rightarrow \lambda_1 \geq \lambda_2$'' is easy to obtain
\emph{without assuming {\bf (U1)}}, just as a straightforward consequence of the comparison arguments.
It is therefore true as soon as $F$ and $L$ satisfy the conditions of the comparison result, i.e.
under far weaker assumptions than the result of Theorem~\ref{unierg}.
 The key point in Theorem~\ref{unierg} is really the result ``$\mu_1 < \mu_2 \Rightarrow \lambda_1 > \lambda_2$''\,.}
\end{remark}

\medskip

Now we turn to the checking of {\bf (U1)} which can be formulated in the following way.

\begin{theor}\label{concunierg} The boundary ergodic cost $\mu$ is unique in the two following cases
\begin{enumerate}
\item[(i)] under the assumption of Theorem~\ref{ergnl},
\item [(ii)] under the assumption of Theorem~\ref{ergl} on $F$ and of Theorem~\ref{ergnl} on $L$, if $F(x,p,M)$ is convex in $(p,M)$ and $L(x,p)$ is convex in $p$.
\end{enumerate}
\end{theor}
\par
\medskip
It is worth mentionning that, in the case of the result (ii), we have the uniqueness of $\mu$ for problems for which we do not have a priori an existence result.

\medskip

\noindent{\bf Proof of Theorem \ref{concunierg}.} In order to apply Theorem~\ref{unierg}, it is enough to check that
 {\bf (U1)} holds.
\par
In the case when {\bf (F1)-(F2)} holds, recalling that we may assume $\mathcal{O} \subset \{x_1>0\}$, we set
$$ w_\e= w - \e \varphi(x)\quad\hbox{for  }x\in \Ob\; ,$$
where $ \varphi(x):= 2-\exp(-\sigma x_1)$ for some $\sigma>0$ choosen later. If $e_1:= (1,0, \cdots, 0)$, denoting by $\ell(x) := \exp(-\sigma x_1)$, we have
\begin{eqnarray*}
F(x,Dw_\e,D^2 w_\e) & = & F(x, Dw-\e \sigma \ell(x)e_1, D^2 w + \e\sigma^2 \ell(x) e_1 \otimes e_1)\\
& \leq & F(x,Dw,D^2 w) -\kappa \e\sigma^2 \ell(x) + K\e \sigma
\ell(x)\; .
\end{eqnarray*}

By choosing $\sigma > K\kappa^{-1}$,  the quantity $-\kappa \e\sigma^2 \ell(x) + K\e \sigma \ell(x)$ becomes strictly negative on $\Ob$ and we have $F(x,Dw_\e,D^2 w_\e)\leq \lambda_\e < \lambda$. The checking for the boundary condition is straightforward using {\bf (L2)}.

In the case when {\bf (F2')} holds, we cannot argue in the same way. We set
$$ w_\e= (1-\e)w - \e c\varphi(x)\; ,$$
where $\varphi(x)$ is defined as above and $c>0$ will be chosen later. By the convexity of $F$, we have
$$ F(x,Dw_\e,D^2 w_\e) \leq (1-\e)F(x,Dw,D^2 w)+\e F(x,-cD\varphi(x), -cD^2 \varphi(x))\; .$$
To conclude, it is enough to show that we can choose $\sigma$ and $c$ in order that
$$F(x,-cD\varphi(x), -cD^2 \varphi(x))<\lambda,~~\mbox {on $\Ob$\,.}
$$
We have
$$
F(x,-cD\varphi(x), -cD^2 \varphi(x)) = F(x,-c\sigma \ell(x)e_1, c
\sigma^2 \ell(x) e_1 \otimes e_1) \; ,$$ and by {\bf (F2')}
\begin{eqnarray*}
F(x,-c\sigma \ell(x)e_1, c \sigma^2 \ell(x) e_1 \otimes e_1)
 &\leq&  F(x,- c\sigma \ell(x)e_1, 0)\\
 &&-\kappa c \sigma^2 \ell(x) \left(\langle\widehat{A^{-1}(x)e_1},e_1\rangle\right)^2+c \sigma^2 \ell(x)o(1)\; .
 \end{eqnarray*}

Finally by using {\bf (F1)} and the fact that $A$ is positive
definite we get
\begin{eqnarray*}
 F(x,-cD\varphi(x), -cD^2 \varphi(x))& \leq&
   F(x,0, 0)+ Kc\sigma \ell(x) \\ & &
   -C \kappa c \sigma^2 \ell(x)+c \sigma^2 \ell(x)o(1)\; ,
   \end{eqnarray*}
for some (small) constant $C>0$ and $o(1)\to 0$ as
$c,\sigma\to\infty.$
 We conclude by first choosing $\sigma$ large
enough and then $c$ large enough. The checking for $L$ is done in
an analogous way and even simpler because we do not need a
sign.\sn

\medskip

Now we turn to an almost immediate corollary of the uniqueness

\begin{cor}\label{mu-cont} { Under the assumptions of either Theorem~\ref{ergnl} or Theorem~\ref{ergl} and Theorem~\ref{concunierg}(iii)}, the map $\lambda \mapsto \mu(\lambda)$ is continuous and decreasing.
\end{cor}
\noindent{\bf Proof.} The solutions $u:=u(\lambda)$ of (\ref{eqn})-({\ref{bc}) we build in the proofs of Theorems \ref{ergnl} and \ref{ergl} with the property $u(x_0) = 0$ are bounded in $C^{0,\beta}(\Ob)$ for $\lambda$ bounded. By Ascoli's Theorem, this means that the $u(\lambda)$ are in a compact subset of $C(\Ob)$ if $\lambda$ remains bounded. Using this property together with the stability result for viscosity solutions and the fact that $\mu$ is also bounded if  $\lambda$ is bounded by the basic estimates on $\alpha u$ of the existence proof, yields easily the continuity of $\mu$ w.r.t $\lambda$. Here, of course, the uniqueness property for $\mu$ plays a central role.

The monotonicity is a direct consequence of Theorem~\ref{unierg}
since it shows that if $\lambda_1 \leq \lambda_2$, then
necessarily $\mu(\lambda_1)\geq \mu(\lambda_2)$. Thus the result
follows.\sn

\begin{cor}\label{lamu} Under the assumptions of Corollary~\ref{mu-cont}, there exists a unique $\lambda:=\tilde \lambda$ such that $\mu(\tilde \lambda)=\tilde \lambda$.
\end{cor}

\noindent{\bf Proof.} The map $\chi(\lambda) :=\lambda - \mu(\lambda)$ is continuous, strictly increasing on $\r$ and satisfies $\chi(-\infty)=-\infty$ and $\chi(+\infty)=+\infty$. Hence the result is a direct consequence of the Intermediate Values Theorem.\sn

\par

\vskip0.5cm

We conclude this section by a result describing a little bit more precisely the dependence of $\mu$ in $F$ and $L$. Of course, since $\lambda$ can be incorporated in $F$, this result gives also informations on the behavior of $\mu$ with respect to $\lambda$ but we argue here with a fixed $\lambda$. We use the natural notation $\mu(F, L)$ to emphasize the dependence of $\mu$ in these two variables.

\begin{theor} If $F_1, F_2$ and $L_1,L_2$ satisfies the assumptions of Corollary~\ref{mu-cont} and if $F_1-F_2$, $L_1-L_2$ are bounded, there exists a constant $\tilde C>0$ such that
$$ |\mu(F_1, L_1) - \mu(F_2,L_2) | \leq \tilde C\left( ||ÊF_1-F_2||_\infty + || L_1-L_2||_\infty\right) \; .$$
\end{theor}

\noindent{\bf Proof.} We start by the {\bf uniformly elliptic case}.
\par
We denote by $u_1$ the solution associated to $F_1, L_1$ and $\mu(F_1, L_1)$. Applying readily the computations of the proof of Theorem~\ref{concunierg}, it is easy to show that $w:=u_1-k ||F_1-F_2||_\infty\varphi$  ($\varphi$ being the function defined in the proof of Theorem \ref{unierg}) is a subsolution for the equation $F_2$. Moreover
$$ L_2(x,Dw) \leq \mu(F_1, L_1) + || L_1-L_2||_\infty + C||F_1-F_2||_\infty\; ,$$
for some constant $C$.
Applying Theorem~\ref{unierg}, we deduce that
$$  \mu(F_2, L_2) \geq  \mu(F_1, L_1) + || L_1-L_2||_\infty + C||F_1-F_2||_\infty\; ,$$
and the result follows by exchanging the roles of $(F_1, L_1)$ and $(F_2, L_2)$.

For the {\bf convex, non uniformly elliptic case}, we argue similarly but by taking this time $w:=\theta u_1-(1-\theta)k\varphi$ with $k>0$ large to be chosen later and for some suitable $0<\theta<1$.  Because of the convexity of $F_2$, $w$ satisfies for some $\lambda>0, C>0$
$$ F_2 (x,Dw,D^2w) \leq ||F_1-F_2||_\infty + \theta \lambda -(1-\theta)Ck\; .$$
We choose $k>0$ large enough and then $\theta$ in order to have
$$ ||F_1-F_2||_\infty + \theta \lambda -(1-\theta)Ck = \lambda\; .$$
Next we examine the boundary condition~:  using again the convexity of $L_1$, we obtain
$$ L_2(x,Dw) \leq \theta \mu(F_1, L_1) + (1-\theta){\overline C}k + ||L_1 - L_2||_\infty\; .$$
As above we deduce
$$ \mu(F_2, L_2) \geq \theta \mu(F_1, L_1) + (1-\theta){\overline C}k + ||L_1 - L_2||_\infty\; .$$
In order to conclude, we have to play with $k$ and $\theta$. The above inequality can be rewritten as
$$ \mu(F_2, L_2) -\mu(F_1, L_1)  - ||ÊL_1 - L_2||_\infty \geq(1-\theta)\left[ {\overline C}k -  \mu(F_1, L_1)\right] \; ,$$
and with the choice of $k$ and $\theta$
$$  (1-\theta)= \frac{ ||ÊF_1-F_2||_\infty}{Ck + \lambda}\; .$$
Finally
$$ \mu(F_2, L_2) -\mu(F_1, L_1)  - ||ÊL_1 - L_2||_\infty \geq ||F_1-F_2||_\infty \frac{{\overline C}k -  \mu(F_1, L_1) }{Ck + \lambda} \; ,$$
and the conclusion follows by letting $k$ to $+ \infty$.~\hfill\sn

\medskip
We conclude this section by showing that, under the hypotheses of Theorem~\ref{ergnl}, the solution of (\ref{eqn})-(\ref{bc}) is unique up to additive constants.
\begin{theor}\label{uniu}
Under the assumptions of Theorem~\ref{ergnl}, the solution of the problem (\ref{eqn})-(\ref{bc}) is unique up to additive constants.
\end{theor}

\medskip
{\noindent\bf Proof.} Suppose by contradiction that $u_1$ and $u_2$ are two solutions of (\ref{eqn})-(\ref{bc}) associated to $\lambda$ and $\mu(\lambda)$, such that the function $w:=u_1-u_2$ is not constant.

We first show that $w$ is a subsolution of a suitable Neumann
problem; this is the aim of the following lemma in which, for
$x\in\dO$, we denote by $D_T w (x)$ the quantity $Dw(x) -
(Dw(x)\cdot n(x))n(x)$. $D_T w (x)$ represents the projection of
$Dw(x)$ on the tangent hyperplane to $\dO$ at $x$. For $X \in
{\mathcal{S}}^n$, we use also the notation
$$ \mathcal{M}^+(X) = \sup_{\kappa Id\leq A \leq K Id}\, {\rm Tr}(AX)\; ,$$
for the Pucci's extremal operator associated to the constants $K$ and $\kappa$ appearing in assumptions {\bf (F1)} and {\bf (F2)} respectively.

\begin{lemma}\label{linea}
Under the assumptions of Theorem~\ref{ergnl},  $w=u_1-u_2$ is a viscosity subsolution
of
\begin{eqnarray}
  -\mathcal{M}^+(D^2 w)-K|Dw|=0&~~ \mbox{in $\O$}\label{eqlinz}\\
\frac{\partial w}{\partial n}-C|D_T w|=0 & ~~\mbox{on $\dO$}\label{bclinz}
\end{eqnarray}
where $C>\max{(K,\displaystyle \frac{\bar K}{\nu})}$, $K,\bar K,\nu$ being the constants appearing in
{\bf (F1)} and {\bf (L1)-(L2)}\,.
\end{lemma}

\medskip
We postpone the (sketch of the) proof of this lemma to the Appendix and conclude the proof of Theorem~\ref{uniu}.
Using this lemma, the function $w=u_1-u_2$ is a non-constant viscosity subsolution of \rec{eqlinz}-\rec{bclinz}. To obtain the contradiction, we use the same arguments as in the step 3 of the proof of Theorem~\ref{ergnl}~: by the Strong Maximum Principle, $w$ cannot achieve its maximum in $\O$. But then Lemma~\ref{SMPbc} and the same arguments as in this step 3 leads to a contradiction.~\hfill\sn

\section{Asymptotic behavior as $t\to +\infty$ of solution of nonlinear equations}\label{ab}
We describe in this section two properties related on the asymptotic behavior of solutions of parabolic equations which are connected to the boundary ergodic cost.

We first consider the evolution problem
\begin{eqnarray}
\chi_t + F(x,D\chi,D^2 \chi)=\lambda & \mbox{in $\O\times (0,\infty)$,}\label{evol1} \\
L(x,D\chi)= \mu & \mbox{on $\partial{\mathcal{O}} \times (0,\infty)$,} \label{evol2} \\
 \chi(x,0) = u_0 (x) & \mbox{in ${\mathcal{O}}$\,.} \label{evol3} \end{eqnarray}

\begin{theor}\label{T:evol1} Under the assumptions of Corollary~\ref{mu-cont}, there exists a unique viscosity solution $\chi$ of (\ref{evol1})-(\ref{evol2})-(\ref{evol3}) which is defined for all time.
Moreover,  $\chi$ remains uniformly bounded in time if and only if $\mu = \mu(\lambda)$.
\end{theor}

\noindent{\bf Proof.} The existence and uniqueness of $\chi$ is a standard result. Only the second part of the result is new.
To prove it, we first assume that $\mu = \mu(\lambda)$. If $u$ is the solution of (\ref{eqn})-(\ref{bc}), it is also a solution of (\ref{evol1})-(\ref{evol2})-(\ref{evol3}) with initial data $u$ and by standard comparison argument
$$ ||Ê\chi(\cdot,t) - u(\cdot)||_\infty \leq ||u_0 - u||_\infty\; ,$$
which implies the claim.

Conversely, if $\chi$ is uniformly bounded, by considering the functions
$$ \chi_\alpha (x,t) :=\chi (x,\alpha^{-1}t)\; ,$$
for $\alpha >0$ small,
it is straigntforward to show that
$$ \overline \chi := {\limsup}^* \chi_\alpha \; \hbox {  and  }\;\underline \chi := {\liminf}_* \chi_\alpha \; ,$$
are respectively sub and supersolution of (\ref{eqn})-(\ref{bc}). A simple application of Theorem~\ref{unierg} shows that $\mu = \mu(\lambda)$. And the proof is complete.~\hfill\sn

\medskip

We next consider the problem
\begin{eqnarray}
\phi_t + F(x,D\phi,D^2 \phi)=0 & \mbox{in ${\mathcal{O}} \times (0,\infty)$,} \label{eqn-evol} \\
\phi_t + L(x,D\phi)= 0 & \mbox{on $\dO \times (0,\infty)$,} \label{bc-evol}\\
 \phi(x,0) = \phi_0 (x) & \mbox{in ${\mathcal{O}}$,} \label {id-evol}
 \end{eqnarray}
where $\phi_0\in C(\Ob)$.

Our result is the
\begin{theor}\label{T:evol2} Under the assumptions of Corollary~\ref{mu-cont}, there exists a unique viscosity solution of
(\ref{eqn-evol})-(\ref{bc-evol})-(\ref{id-evol}) which is defined for all time. Moreover, as $t \to +\infty$, we have
$$ \frac{\phi(x,t)}{t} \to - \tilde \lambda\quad\hbox{uniformly on  }\Ob\; ,$$ where $\tilde \lambda$ is defined in Corollary~\ref{lamu}
\end{theor}

\noindent{\bf Proof.} We denote by $\tilde u$ the solution of (\ref{eqn})-({\ref{bc})  associated to $\lambda = \tilde \lambda$ and $\mu = \tilde \lambda$.

The existence and uniqueness of $\phi$ is a consequence of the results in \cite{b2}. Moreover, since $\tilde u-\tilde \lambda t$ is a solution of (\ref{eqn-evol})-(\ref{bc-evol}), the comparison result for this evolution equation yields
$$ || \phi(x,t) - \tilde u(x)+\tilde \lambda t||_\infty \leq ||\phi_0 - \tilde u ||_\infty\; .$$
Dividing by $t$ and letting $t$ tends to infinity provides the result.\sn

\section{On ergodic stochastic control problems}\label{escp}

We are interested in this section in control problems of diffusion processes with reflection. The dynamic is given by the solution of the following problem in which the unknown is a pair
$((X_t)_{t\geq 0} ,(k_t)_{t\geq 0})$ where $(X_t)_{t\geq 0}$ is a
continuous process
in $\R$ and $(k_t)_{t\geq 0}$ is a process with bounded variations
\begin{equation}\label{ROe}
\left\{ \begin{array}{ll}
dX_t = b(X_t, \alpha_t)dt + \sigma (X_t,\alpha_t)dW_t
- dk_t\; , & X_0=x \in \Ob , \\  k_t =\int_0^t 1_{\dO} (X_s)
\gamma (X_s)d|k|_s \; ,& X_t \in \Ob \; ,\quad \forall t\geq 0
\; , \\
\end{array} \right.
\end{equation}
where $(W_t)_t$ is a $p$-dimensional Brownian motion for some $p\in \N$. The process
$(\alpha_t)_t$, the {\em control}, is some progressively  measurable process with respect to the filtration associated to the Brownian motion with values in  a compact metric space $\mathcal{A}$.  The drift $b$ and the diffusion
matrix $\sigma$ are continuous functions defined on $\overline \Omega
\times \mathcal{A} $ taking values respectively in $\R$ and in the
space of $N \times p$ matrices. We assume that both $b$ and $\sigma$ are Lipschitz continuous in $x$, uniformly in $\alpha \in \mathcal{A}$. Finally $\gamma$ satisfies the assumptions given in Section~\ref{S:lbc}.

Under these assumptions, there exists a unique pair $((X_t)_{t\geq 0} ,(k_t)_{t\geq 0})$  solution of this problem, the existence being proved in Lions \& Sznitman \cite{LiSz} and the uniqueness in Barles \& Lions \cite{BaLi}.

Then we define the value-function of the finite horizon, stochastic control problem by
\begin{equation}\label{valuefunct}
U (x,t) = \inf_{\control}\,\Ex\left[
\int_0^t [f(X_s,\alpha_s)+\lambda ] dt + \int_0^t [g(X_s)+\mu]d|k|_s + u_0 (X_t) \,\right] \; ,
\end{equation}
where $\Ex$ denotes the conditional expectation with respect to the
event $\{X_0=x\}$, $f$ is a continuous function defined on $\Ob\times \mathcal{A}$ which is Lipschitz continuous in $x$ uniformly w.r.t $\alpha \in \mathcal{A}$, $g \in C^{0,\beta}(\dO)$ and $u_0\in C(\Ob)$, $\lambda$ and $\mu$ are constants.
\par
Under the above assumptions, by classical results, $U$ is the unique viscosity solution of
\begin{eqnarray*}
U_t + F(x,DU,D^2 U)=\lambda & \mbox{in $\O\times (0,\infty)$,} \\
\frac{\partial U}{\partial \gamma}= g+\mu & \mbox{on $\partial{\mathcal{O}} \times (0,\infty)$,} \\
U(x,0) = u_0 (x) & \mbox{in ${\mathcal{O}}$,} \end{eqnarray*}
with
$$ F(x,p,M)=\sup_{\alpha\in \mathcal{A}}\,
\left\{-\frac{1}{2}{\rm Tr}[
a(x, \alpha)M]- \langle b(x,\alpha), p\rangle -f(x,\alpha)\right\}
$$
for any $x\in \Ob$, $p\in \R$ and $M\in {\mathcal{S}}^n$ where $a(x, \alpha)= \sigma(x,\alpha)\sigma^T(x, \alpha)$.
We are going to use this Hamilton-Jacobi-Bellman type evolution problem both to study the stationary ergodic problem (and, in particular, to revisit the result of Theorem~\ref{unierg} in a \textit{degenerate} context) and to connect the constant $\mu (\lambda)$ with the behavior of $U$ as $t\to \infty$ in the spirit of Theorem~\ref{T:evol1}.
Our result is the following.
\par\medskip
\begin{theor}\label{T:condk} Under the above assumptions on $\sigma$, $b$, $f$, $g$ and $u_0$, we have
\begin{enumerate}
  \item[(i)] For the stationary ergodic problem, the analogue of Theorem~Ê\ref{unierg} (i.e. ``$\mu_1 < \mu_2 \Rightarrow \lambda_1 > \lambda_2$'') is equivalent to the property
\begin{equation}\label{condk}
\sup_{x\in \Ob} \limsup_{t\to +\infty} \left( \inf_{\control}\Ex \int_0^{t} d|k|_s \right) = +\infty\; .
\end{equation}
In particular, under this condition, if $\mu(\lambda)$ exists for some $\lambda \in \r$, it is unique.
 \item[(ii)] If (\ref{condk}) holds, for any $\lambda$, there exists at most a constant $\mu(\lambda)$ for which $U$ is uniformly bounded.
\item[(iii)] We set
\begin{equation}\label{m}
m(x,t) := \inf_{\control}\Ex \left( \int_0^{t} d|k|_s \right)\; .
\end{equation}
Assume that (\ref{condk}) holds and that there exists a constant $\mu(\lambda)$ for which $U$ is uniformly bounded.
If $m(x_n,t_n)\to + \infty$ with $x_n \in \Ob$, $t_n\to +\infty$, then
\end{enumerate}
\begin{equation}\label{formu}
\mu(\lambda):= -  \lim_{n \to +\infty}\left\{ \inf_{\control}\,\left[\left({I\!\!E_{x_n}} \left[\int_0^{t_n} d|k|_s\right]\right)^{-1}J(x_n,t_n, \control)\right]\right\}.
\end{equation}
where
$$J(x_n,t_n, \control):= {I\!\!E_{x_n}}
 \left(
\int_0^{t_n} [f(X_s,\alpha_s)+\lambda ] dt + \int_0^{t_n} g(X_s)d|k|_s \right)$$
\end{theor}

This result gives a complete characterization of $\mu(\lambda)$ when it exists and it points out the conditions under which this constant is unique. In particular, (\ref{condk}) is a justification of the idea that in order to have a unique $\mu(\lambda)$, the boundary condition has to be ``sufficiently seen''.

Of course, the weak part of this result is the existence of $\mu(\lambda)$: unfortunately, in this case, we cannot have a better result than the uniformly elliptic case since {\bf (F2')} leads to assume that the equation is uniformly elliptic. Indeed, if one considers {\bf (F2')} with $N = cq\otimes q$ where $q=\widehat {A^{-1}(x)p}$ and $c>0$ is very large, then by dividing by $c$ and letting $c$ tends to $+ \infty$, we are lead to
$$ \sup_{\alpha\in \mathcal{A}}\,
\left [-\frac{1}{2}\langle a(x, \alpha)q,q\rangle \right]\leq - \kappa \; ,$$
since $|q| = 1$. In other words, for any $x \in \Ob$ and $\alpha\in \mathcal{A}$,
$\langle  a(x, \alpha)q,q\rangle  \geq \kappa$. Since this has to be true for any $p$, hence for any $q$, this shows that the equation has to be uniformly elliptic.

This uniform elliptic case is the purpose of the following corollary.

\begin{cor}\label{conclsto}Under the above assumptions on $\sigma$, $b$, $f$, $g$ and $u_0$ and if there exists $\nu >0$ such that $a(x,\alpha) \geq \nu Id$ for any $x\in \Ob$ and $\alpha \in \mathcal{A}$, then (\ref{condk}) holds and for any $\lambda \in \r$, there exists a unique $\mu(\lambda)\in \r$ for which $U$ is uniformly bounded. This constant $\mu(\lambda)$ is given by (\ref{formu}) and it is also the unique constant for which the associated stationary Bellman boundary value problem has a solution.
\end{cor}

We skip the proof of this result since it follows easily from either Theorem~\ref{ergnl} or \ref{ergl}, Theorem~\ref{unierg} and \ref{concunierg} and Theorem~\ref{T:condk}.

Before turning to the proof of Theorem~\ref{T:condk}, we want to point out that, in general, even if (\ref{condk}) holds, $m$ is not expected to converge to infinity uniformly on $\Ob$, nor even at any point of $\Ob$. Indeed it is very easy, in particular in the deterministic case, to build situations for which the drift is like $n$ in a neighborhood of $\dO$ (and therefore the trajectory are pushed to $\dO$ leading to (\ref{condk})) while $b$ can be identically $0$ inside $\mathcal{O}$ and therefore for such points $k_s \equiv 0$. As a consequence of this remark, the admittedly strange formulation of (iii) cannot be improved.

\medskip

\noindent{\bf Proof.} We first prove (i). We first assume that the
analogue of Theorem \ref{unierg} holds and we want to show that
(\ref{condk}) holds. We argue by contradiction assuming that it
does not; this implies that the function $m$ defined in \rec{m} is
uniformly bounded. Indeed (\ref{condk}) is clearly equivalent to
$$ \sup_{x\in \Ob} \limsup_{t\to +\infty} \, m(x,t) = +\infty\; ,$$
and the function $m$ is increasing in $t$.

We choose above $f=g=u_0 = 0$. Arguing as in the proof of Theorem~\ref{T:evol1}, we see that, as $t\to \infty$, $\underline m := {\liminf}_*\,m$ is a supersolution of the stationary equation with $\lambda = 0$ and $\mu =1$ while $0$ is a solution of this problem with $\lambda = 0$ and $\mu =0$. This is a contradiction with the assumption.

Conversely, if (\ref{condk}) holds, let $u_1$ be an usc subsolution of the stationary problem associated to $\lambda_1, \mu_1$ and $u_2$ be a lsc supersolution of the stationary problem associated to $\lambda_2, \mu_2$, with $\mu_1 < \mu_2$. The functions $u_1$ and $u_2$ are respectively sub and supersolution of the evolution equation (with, say, initial datas $||u_1||_\infty$ and $-||u_2||_\infty$ respectively); therefore, for any $x\in \Ob$ and $t>0$
$$u_1 (x) \leq \inf_{\control}\,\Ex\left[
\int_0^t [f(X_s,\alpha_s)+\lambda_1 ] dt + \int_0^t [g(X_s)+\mu_1]d|k|_s + ||u_1||_\infty \,\right] \; ,$$

$$u_2 (x) \geq \inf_{\control}\,\Ex\left[
\int_0^t [f(X_s,\alpha_s)+\lambda_2] dt + \int_0^t [g(X_s)+\mu_2]d|k|_s - ||u_2||_\infty \,\right] \; .$$

\noindent Let us take a sequence $(x_n,t_n)\in\overline\O\times(0,+\infty)$ such that
$t_n\to+\infty$ and $m(x_n,t_n)\to +\infty$ as $n\to +\infty.$  Let $\alpha_n$ be an $\varepsilon$-optimal control for the ``inf'' in the $u_2$ inequality with $\varepsilon =1$. Using also $\alpha_n$ for $u_1$ and subtracting the two inequalities, we obtain
$$ (u_1-u_2)(x_n) \leq {I\!\! E_{x_n}}\int_0^{t_n} [\lambda_1 -\lambda_2] dt + \int_0^{t_n} [\mu_1-\mu_2]d|k|_s + O(1)\; .$$
In this inequality, by the definition of $m$, the $k$-term is
going to $-\infty$ since $\mu_1-\mu_2 < 0$ but the left-hand side
is bounded; so necessarily $\lambda_1 -\lambda_2>0$.\par
 We next prove (ii). Suppose by contradiction that there are $\mu_1$ and $\mu_2$ such that
the corresponding value functions $U_1$ and $U_2$ defined by (\ref{valuefunct}) are uniformly
bounded in $\overline\O\times [0,\infty).$  We assume that $\mu_1>\mu_2$. We have
\begin{equation}
U_1(x,t)-U_2(x,t)\ge (\mu_1-\mu_2) \inf_{\control} \left({I\!\!E_{x}} \int_0^{t} d|k|_s\right)\,.
\end{equation}
By letting $t\to  +\infty$ we get a contradiction because of the condition \rec{condk}.

We leave the proof of   (iii) to the reader since it is an easy adaptation of the arguments we give above.
~\hfill\sn

\section{Appendix}

\subsection{A comparison argument using only {\bf (F1)-(F2)}}\label{ca}

The difficulty comes from {\bf (F1)} and can be seen on a term like $-{\rm Tr}(A(x)D^2 u)$: in general, one assumes that $A$ has the form $A=\sigma\sigma^T$ for some Lipschitz continuous matrix $\sigma$ and the uniqueness proof uses $\sigma$ in an essential way, both in the degenerate and nondegenerate case. Here we want just to assume $A$ to be nondegenerate and Lipschitz continuous and we do not want to use $\sigma$, even if, in this case, the existence of such $\sigma$ is well-known.

In the comparison argument of \cite{b2}, the only difference is in the estimate of the difference $F(x,p,X)-F(y,q,Y)$.

The key lemma in \cite{GBMR} to solve this difficulty is the following: if the matrices $X,Y$ satisfy (\ref{inegmatthbis}) (with $\eta = 0$) then
$$ X-Y \leq - \frac{\tilde K \e^2}{6} (t X + (1-t)Y)^2 \quad \hbox{for all  } t \in[0,1]. $$
A slight modification of this argument allows to take in account the $\eta$ term and yields
$$ X-Y \leq - \frac{\tilde K \e^2}{6} (t X + (1-t)Y)^2 + O(\eta) \quad \hbox{for all  } t\in [0,1]. $$
Now we show how to estimate $F(x,p,X)-F(y,q,Y)$. By using {\bf (F1)-(F2)} together with the above inequality for $t=0$, we get
\begin{eqnarray*}
F(x,p,X)-F(y,q,Y)& \geq &F(x,p,X) - F(x,p,Y + O(\eta)) \\
& & - K(|p-q| + |x-y|(|p|+||Y||) + O(\eta)   \\
 & \ge & \kappa  \frac{\tilde K \e^2}{6} {\rm Tr}(Y^2)- K(|p-q| + |x-y|(|p|+||Y||) + O(\eta)\; .
\end{eqnarray*}
In this inequality, the ``bad'' term is $ K|x-y|||Y||$ since the
estimates on the test-function does not ensure that it converges
to $0$. But this term is controlled by the ``good term''  ${\rm
Tr}(Y^2)$ in the following way: by Cauchy-Schwarz's inequality
$$  K |x-y|||Y|| \geq -\kappa  \frac{\tilde K \e^2}{6} {\rm Tr}(Y^2) - O\left(\frac{|x-y|^2}{\e^2}\right)\; .$$
And this estimate is now sufficient since we know that $\displaystyle \frac{|x-y|^2}{\e^2}\to 0$ as $\e \to 0$.

\subsection{Proof of Lemma~\ref{SMPbc}}

We use here argument which are borrowed from \cite{BaDL}. We prove the result under the weaker assumption {\bf (F1)-(F2')}.
\par
Since ${\mathcal{O}}$ is a $C^2$ domain, for $s>0$ small enough, $d(x-sn(x))=s$ where $d$ is the distance to the boundary $\dO$.
We set $x_0 = x-sn(x)$ for such an $s$ and we build a function $\varphi$ of the following form
$$\varphi(y) = \exp(- \rho s^2)-\exp(-\rho |y-x_0|^2)\; ,$$
where $\rho$ has to be chosen later.
Finally we choose $r=s/2$.
Since $s=|x-x_0|$, we have $\varphi(x) =0$ and if $y\in \dO \cap \overline B(x,r)-\{x\}$, $|y-x_0| \ge s/2$ and therefore $\varphi(y)>0$. Moreover
$$ D\varphi(y) = 2\rho(y-x_0)\exp(-\rho |y-x_0|^2)\; ,$$
and by the definition of $x_0$, $D\varphi(x) = kn(x)$ with $k=2s\rho\exp(-\rho s^2)>0$.
Finally, we compute  $F_\infty (y,D \varphi (y), D^2 \varphi(y)).$
Using the notations $\ell(y) =  2\rho\exp(-\rho |y-x_0|^2)$ and $p(y) =y-x_0$, we have
$$ F_\infty (y,D \varphi (y), D^2 \varphi(y))=
F_\infty (y,\ell(y)p(y), \ell(y) Id - 2\rho \ell(y)  p(y) \otimes p(y))\; .$$

By homogeneity, it is enough to have
$$ F_\infty (y, p(y),  Id - 2\rho p(y) \otimes p(y))>0\; .$$
We notice that, in $B(x,r)$, $p(y)$ does not vanish and {\bf
(F2')} yields
\begin{eqnarray*}
 F_\infty (y, p(y),  Id - 2\rho p(y) \otimes
p(y))&\geq&
 2\kappa\rho\langle \widehat{A^{-1}(y)p(y)},p(y)\rangle
^2 \\ && + F_\infty (y, p(y),  Id) +o(1)2\rho|p(y)|^2\; .
\end{eqnarray*}
In order to have the left-hand side positive,
it is enough to choose $\rho$ large enough. And the proof is complete.~\hfill\sn

\subsection{Sketch of the Proof of Lemma~\ref{linea}}

We just sketch the proof since we follow very closely the strategy of proof of Lemma 2.6 in Arisawa \cite{A3}.  Let $\phi\in C^2(\overline\O)$ be such that $w-\phi$ has a local maximum at $\bar x\in\overline\O.$ We suppose that $\bar x\in\dO$, the case $x\in\O$ being similar and even simpler.

For all $\eps>0$ and $\eta>0$, we introduce the auxiliary function
\begin{equation}\label{aux}
\Phi_{\eps,\eta}(x,y)=u_1(x)-u_2(x)-\psi_{\eps,\eta}(x,y)-\phi(\frac{x+y}{2})-|x-\bar x|^4
\end{equation}
where $\psi_{\eps,\eta}(x,y)$ is the test function built in Barles \cite{b2} relative to the boundary condition \rec{bc}.
Let $(x_\eps, y_\eps)$ be the maximum point of $\Phi_{\eps,\eta}(x,y)$ in $\overline\O\times\overline\O$. Since $\bar x$ is a strict local maximum point of $x\mapsto w(x) -\phi(x)-|x-\bar x|^4$, standard arguments show that
$$(x_\eps,y_\eps)\to (\bar x,\bar x)\,\hbox{  and  }\,\frac{|x_\eps-y_\eps|^2}{\eps^2}\to 0~~ \mbox{as $\eps\to 0.$}
$$
On the other hand, by construction we have
$$ L(x_\eps, D_x\psi_{\eps,\eta}(x_\eps,y_\eps))>\mu~~\mbox{if $x_\eps\in\partial\O\,,$}$$
$$ L(y_\eps, -D_y\psi_{\eps,\eta}(x_\eps,y_\eps))<\mu~~\mbox{if $y_\eps\in\partial\O\,.$}$$
Moreover, if $\displaystyle \zeta_{\eps,\eta}(x,y):=\psi_{\eps,\eta}(x,y)+\phi(\frac{x+y}{2})+|x-\bar x|^4$, by standard arguments (cf. \cite{cil}), we know that, for every $\alpha>0$, there exist $X,Y\in{\mathcal{S}^n}$ such that
$$
(D_x\zeta_{\eps,\eta}(x_\eps,y_\eps),X)\in \overline{J}^{2,+}_{\overline\O}u_1(x_\eps)\,,
$$
$$
(-D_y\zeta_{\eps,\eta}(x_\eps,y_\eps),Y)\in \overline{J}^{2,-}_{\overline\O}u_2(y_\eps)\,,
$$
and
$$
-(\frac{1}{\alpha}+||D^2 \zeta_{\eps,\eta}(x_\eps,y_\eps)||)Id\leq
\left(\begin{array}{cc} X & 0 \\ 0 & -Y \end{array}\right)
  \leq (Id+{\alpha}
D^2\zeta_{\eps,\eta}(x_\eps,y_\eps))D^2\zeta_{\eps,\eta}(x_\eps,y_\eps)\, .
$$

Now suppose that
$$ \frac{\partial\phi }{\partial n}(\bar x)-C|D_T \phi(\bar x)|>0\, .$$
If $x_\eps \in\dO$, then, for $\eps$ small enough, we have
$$ L(x_\eps,D_x \zeta_{\eps,\eta}(x_\eps,y_\eps))\ge L(x_\eps, D_x \psi_{\eps,\eta})
+ \frac{1}{2}(\nu \frac{\partial \phi(\bar x)}{\partial n}-K|D_T \phi(\bar x)|)
+o_\eps(1)>\mu\,,$$
while if $y_\eps \in \dO$
$$ L(y_\eps,-D_y \zeta_{\eps,\eta}(x_\eps,y_\eps))\le L(y_\eps, -D_y \psi_{\eps,\eta})
- \frac{1}{2}(\nu \frac{\partial \phi(\bar x)}{\partial n}-K|D_T \phi(\bar x)|)
+o_\eps(1)<\mu\,.
$$
Therefore, if $\eps$ is small enough, wherever $x_\eps,y_\eps$ lie we have
$$
F(x_\eps,D_x \zeta_{\eps,\eta}(x_\eps,y_\eps),X)\le \lambda\,,
$$
$$
F(y_\eps,-D_y \zeta_{\eps,\eta}(x_\eps,y_\eps),Y)\ge \lambda\,,
$$
By subtracting the above  inequalities, using the above estimates
on $X,Y$ together with the arguments of Subsection~\ref{ca}, the
assumption {\bf (F1)} and {\bf (F2)} and the definition of the
Pucci's extremal operator $\mathcal{M}^+$, by letting $\eps$ tend
to $0$, we are lead to
$$ -\mathcal{M}^+(D^2 \phi(\bar x))-K|D\phi(\bar x) | \leq 0\; ,$$
and the conclusion follows.~\hfill\sn

\subsection{The $C^{0,\alpha}$ regularity results and estimates of \cite{BDL}}

As mentioned in the introduction, we describe in this section the results of \cite{BDL} we are using in this paper, in order to have an as self-contained article as possible. In fact, since we use here global estimates (and not local ones), we can follow the remark at the end of the second section in \cite{BDL} and have results with a little bit weaker assumptions. Of course, we reformulate the results of \cite{BDL} in this global framework.

These results concern nonlinear Neumann boundary
value problems of the form
\begin{equation}\label{np} \left\{\begin{array}{ll}
  F(x,u,Du,D^2u)=0 & \mbox{in $O$,} \\
\displaystyle G(x,u,Du)=0& \mbox{on $\partial O $,} \\
   \end{array} \right.
   \end{equation}
where $O\subset\R$ is a smooth, bounded domain, $F$ and $G$ are, at least, real-valued continuous functions defined respectively  on $\overline O \times \r \times\R\times {\mathcal{S}}^n$ and $\partial O\times\r\times\r^n$.

The assumptions are the following~: on the domain, we require
\begin{enumerate}
\item[{\bf (H1)}]~~{\bf (Regularity of the boundary)} $O$ is a bounded
domain with a $C^2$--boundary.\par
\end{enumerate}
while the basic assumptions on $F$ and $G$ are the
\begin{enumerate}
\item[{\bf (H2)}]~~{\bf (Growth Condition on $F$)} For any
$R>0$, there exist positive constants $C_1^R$, $C_2^R$, $C_3^R$
and functions $\omega_1^R$, $\omega_2^R: \r^+ \to \r$ such that
$\omega_1^R (0+)=0$ and $\omega_2^R(r)=O(r)$ as $r\to 0,$  and for any
$x,y \in \overline O$, $-R\leq u,v \leq R$,
$p,q \in \R$, $M \in {\cal{S}}^n$ and $K>0$ \begin{eqnarray*}
F(x,u,p,M)-F(y,v,q,M + KId)     & \leq &  \omega_1^R(|x-y|(1+
|p|+|q|) +|p-q|)||M|| \\ & & + \omega^R_2(K) +
C_1^R  + C_2^R(|p|^2+|q|^2) \\
      &  & +C_3^R|x-y|(|p|^3+|q|^3)\;  .
\end{eqnarray*}
\end{enumerate}
and
\begin{enumerate}
\item[{\bf(G1)}] For all $R>0$, there exists $\mu_{R}>0$ such
that, for every $(x,u,p)\in \partial O\times[-R,R]\times\r^n,$
and $\lambda>0$, we have
\begin{equation}\label{croissbis}
G(x,u,p+\lambda n(x))-G(x,u,p)\ge \mu_{R}\lambda\,,
\end{equation}
\hspace{-\leftmargin}where $n(x)$  denotes the unit outward normal vector to $\partial
O$ at $x \in \partial O$.
\item[{\bf(G2)}] For all $R>0$ there
is  a constant $K_{R} >0$ such that, for all $x,y\in
\partial O,$ $p,q\in\r^n$, $u,v\in [-R,R] $, we have
\begin{equation}
|G(x,u,p)-G(y,v,q)|\le K_R\left[(1+|p|+|q|)|x-y|+|p-q|+|u-v|
\right]\,.
\end{equation}
\end{enumerate}
Before formulating additional assumptions, we want to point out that {\bf (H2)} is obviously satisfied when {\bf (F1)} holds.

These basic assumptions have to be complemented by a ``strong ellipticity assumption'' which are different in the {\em linear case}, i.e. when $G$ is an affine function of $p$, typically when it is of the form
\begin{equation}\label{oblder}
\langle Du,\gamma(x)\rangle +a(x) u(x)+g(x)=0\quad\mbox{on $\partial O$}.
\end{equation}
and in the {\em nonlinear case}.

In the linear case, the ``strong ellipticity assumption'' is the following
\begin{enumerate}
\item[{\bf (H3a)}] {\bf  Oblique-derivative boundary
condition and ellipticity~:} there exists a \linebreak Lipschitz continuous function $A\colon\overline O\to  {\cal{S}}^n$
with $A\ge c_0Id,$ for some $c_0>0$ such that $A(x)\gamma(x)=n(x)$ for every $x\in
\partial O$, and for any $R>0$, there exist $L_R, \lambda_{R}
> 0$ such that, for all $x\in \overline O,$ $|u|\le R,$
$|p|>L_{R}$ and $M,N\in {\cal{S}}^n$ with $N\geq 0$, we have
\begin{equation}
  \label{ndpsod}
  F(x,u,p,M + N )-F(x,u,p,M) \le -\lambda_{R} \langle N \widehat{A^{-1}(x)p},
\widehat{A^{-1}(x)p}\rangle + o(1)||N|| \, ,
\end{equation}
where $o(1)$ denotes a function of the real variable $|p|$ which converges to $0$ as $|p|$ tends to infinity.
\end{enumerate}

Finally, on the boundary condition (\ref{oblder}), we require
\begin{enumerate}
\item[{\bf (H4)}] {\bf (Regularity of the boundary condition)}
The functions $\gamma$ and $a$ in (\ref{oblder}) are Lipschitz continuous on $\partial O$, $\langle \gamma(x), n(x)\rangle\ge \beta>0$ for any $x\in \partial O$  and  $g$ is
is in $C^{0,\beta}(\partial O)$ for some  $0<\beta \leq 1$.
\end{enumerate}

The result in the {\bf linear case} is the
\begin{theor}\label{theorregul}
Assume {\bf (H1)-(H2)-(H3a)-(H4)}.  Then every continuous
viscosity solution $u$ of \rec{np} with $G$ given by
(\ref{oblder}) is in $C^{0,\alpha}(\overline O)$ for any
$0<\alpha< 1$ if $\beta = 1$ and  with $\alpha=\beta$ if
$\beta<1$. Moreover the $C^{0,\alpha}$--norm of $u$ depend
only on $O$, $F$, $\gamma$, $a$, $g$ through the constants and
functions appearing in {\bf  (H2)-(H3a)}, the
$C^{0,1}$--norm of $\gamma$ and $a$, the $C^{0,\beta}$--norm
of $g$ and the $C^2$--norm of the distance function of the
boundary including the modulus of continuity of $D^2 d$.
\end{theor}

Now we turn to the nonlinear case where we assume uniform ellipticity, namely

\begin{enumerate}
\item[{\bf (H3b)}]~~{\bf (Uniform ellipticity)} For any $R>0$, there is
$\lambda_R > 0$ such that, for all $x\in\overline O$, $-R\leq u \leq
R$, $p \in \R$ and $M,N\in {\cal{S}}^n$ such that $M\le N$, we have $$
F(x,u,p,M)-F(x,u,p,N) \geq \lambda_R\mbox{Tr} (N-M)\;.$$
\end{enumerate}

For the nonlinear boundary condition, we require
\begin{enumerate}
\item[{\bf (G3)}]  For all $R>0$ and $M>0$ there is    $K_{R,M}>0$ such
that
\begin{equation}
|\langle\frac{\partial G}{\partial p}(x,u,p),p\rangle-G(x,u,p)|\le
K_{R,M}\,,
\end{equation}
\hspace{-\leftmargin}for all $x\in\partial O$ and for all $p\in\r^n$, $|p|\ge M,$
$|u|\le R$\,.
   \item[{\bf (G4)}] There is a function $G_\infty\colon\partial
O\times\r\times\r^n\to\r$ such that
\begin{equation}
\frac{1}{\lambda}G(x,u,\lambda p)\to
G_\infty(x,u,p)\,\quad\mbox{as $\lambda\to\infty\,.$}
\end{equation}
\hspace{-\leftmargin} uniformly in $(x,u,p).$
\end{enumerate}

The result in the {\bf nonlinear case} is the
\begin{theor}\label{theorregulbis}
Assume {\bf (H1)-(H2)-(H3b)} and {\bf(G1)-(G4)}.  Then every
bounded continuous solution $u$ of \rec{np} is in
$C^{0,\alpha} (\overline O)$ for any $0<\alpha<1$. Moreover
the $C^{0,\alpha}$--norm of $u$ depend only on $O$, $F$,
$G$,   through the constants and functions appearing in {\bf
(H2)-(H3b)}, and in {\bf(G1)-(G4)}, the $C^2$--norm of the
distance function of the boundary including the modulus of
continuity of $D^2 d$.
\end{theor}

\bigskip

\par\bigskip \centerline{{\sc Acknowledgements}}

\smallskip

  {The second author was partially supported by M.I.U.R., project ``Viscosity, metric, and control theoretic methods for nonlinear partial differential equations'' and by G.N.A.M.P.A, project ``Equazioni alle derivate parziali  e teoria del controllo".

The authors wish to thank the anonymous referee for several valuable suggestions which allow them to improve the form of the article.

\end{document}